\documentclass[12pt]{article} 
\usepackage{latexsym} 
\usepackage{amssymb} 
\usepackage{amsmath} 

\begin{document}


\newtheorem{theorem}{Theorem} 
\newtheorem{problem}{Problem} 
\newtheorem{definition}{Definition} 
\newtheorem{lemma}{Lemma} 
\newtheorem{proposition}{Proposition} 
\newtheorem{corollary}{Corollary} 
\newtheorem{example}{Example} 
\newtheorem{conjecture}{Conjecture} 
\newtheorem{algorithm}{Algorithm} 
\newtheorem{exercise}{Exercise} 
\newtheorem{remarkk}{Remark} 
 
\newcommand{\be}{\begin{equation}} 
\newcommand{\ee}{\end{equation}} 
\newcommand{\bea}{\begin{eqnarray}} 
\newcommand{\eea}{\end{eqnarray}} 
\newcommand{\beq}[1]{\begin{equation}\label{#1}} 
\newcommand{\eeq}{\end{equation}} 
\newcommand{\beqn}[1]{\begin{eqnarray}\label{#1}} 
\newcommand{\eeqn}{\end{eqnarray}} 
\newcommand{\beaa}{\begin{eqnarray*}} 
\newcommand{\eeaa}{\end{eqnarray*}} 
\newcommand{\req}[1]{(\ref{#1})} 
 
\newcommand{\lip}{\langle} 
\newcommand{\rip}{\rangle} 
\newcommand{\uu}{\underline} 
\newcommand{\oo}{\overline} 
\newcommand{\La}{\Lambda} 
\newcommand{\la}{\lambda} 
\newcommand{\eps}{\varepsilon} 
\newcommand{\om}{\omega} 
\newcommand{\Om}{\Omega} 
\newcommand{\ga}{\gamma} 
\newcommand{\rrr}{{\Bigr)}} 
\newcommand{\qqq}{{\Bigl\|}} 
 
\newcommand{\dint}{\displaystyle\int} 
\newcommand{\dsum}{\displaystyle\sum} 
\newcommand{\dfr}{\displaystyle\frac} 
\newcommand{\bige}{\mbox{\Large\it e}} 
\newcommand{\integers}{{\Bbb Z}} 
\newcommand{\rationals}{{\Bbb Q}} 
\newcommand{\reals}{{\rm I\!R}} 
\newcommand{\realsd}{\reals^d} 
\newcommand{\realsn}{\reals^n} 
\newcommand{\NN}{{\rm I\!N}} 
\newcommand{\DD}{{\rm I\!D}} 
\newcommand{\degree}{{\scriptscriptstyle \circ }} 
\newcommand{\dfn}{\stackrel{\triangle}{=}} 
\def\complex{\mathop{\raise .45ex\hbox{${\bf\scriptstyle{|}}$} 
     \kern -0.40em {\rm \textstyle{C}}}\nolimits} 
\def\hilbert{\mathop{\raise .21ex\hbox{$\bigcirc$}}\kern -1.005em {\rm\textstyle{H}}} 
\newcommand{\RAISE}{{\:\raisebox{.6ex}{$\scriptstyle{>}$}\raisebox{-.3ex} 
           {$\scriptstyle{\!\!\!\!\!<}\:$}}} 
 
\newcommand{\hh}{{\:\raisebox{1.8ex}{$\scriptstyle{\degree}$}\raisebox{.0ex} 
           {$\textstyle{\!\!\!\! H}$}}} 

\newcommand{\OO}{\won} 
\newcommand{\calA}{{\cal A}} 
\newcommand{\calB}{{\cal B}} 
\newcommand{\calC}{{\cal C}} 
\newcommand{\calD}{{\cal D}} 
\newcommand{\calE}{{\cal E}} 
\newcommand{\calF}{{\cal F}} 
\newcommand{\calG}{{\cal G}} 
\newcommand{\calH}{{\cal H}} 
\newcommand{\calK}{{\cal K}} 
\newcommand{\calL}{{\cal L}} 
\newcommand{\calM}{{\cal M}} 
\newcommand{\calO}{{\cal O}} 
\newcommand{\calP}{{\cal P}} 
\newcommand{\calU}{{\cal U}} 
\newcommand{\calX}{{\cal X}} 
\newcommand{\calXX}{{\cal X\mbox{\raisebox{.3ex}{$\!\!\!\!\!-$}}}} 
\newcommand{\calXXX}{{\cal X\!\!\!\!\!-}} 
\newcommand{\gi}{{\raisebox{.0ex}{$\scriptscriptstyle{\cal X}$} 
\raisebox{.1ex} {$\scriptstyle{\!\!\!\!-}\:$}}} 
\newcommand{\intsim}{\int_0^1\!\!\!\!\!\!\!\!\!\sim} 
\newcommand{\intsimt}{\int_0^t\!\!\!\!\!\!\!\!\!\sim} 
\newcommand{\pp}{{\partial}} 
\newcommand{\al}{{\alpha}} 
\newcommand{\sB}{{\cal B}} 
\newcommand{\sL}{{\cal L}} 
\newcommand{\sF}{{\cal F}} 
\newcommand{\sE}{{\cal E}} 
\newcommand{\sX}{{\cal X}} 
\newcommand{\R}{{\rm I\!R}} 
\renewcommand{\L}{{\rm I\!L}} 
\newcommand{\vp}{\varphi} 
\newcommand{\N}{{\rm I\!N}} 
\def\ooo{\lip} 
\def\ccc{\rip} 
\newcommand{\ot}{\hat\otimes} 
\newcommand{\rP}{{\Bbb P}} 
\newcommand{\bfcdot}{{\mbox{\boldmath$\cdot$}}} 
 
\renewcommand{\varrho}{{\ell}} 
\newcommand{\dett}{{\textstyle{\det_2}}} 
\newcommand{\sign}{{\mbox{\rm sign}}} 
\newcommand{\TE}{{\rm TE}} 
\newcommand{\TA}{{\rm TA}} 
\newcommand{\E}{{\rm E\,}} 
\newcommand{\won}{{\mbox{\bf 1}}} 
\newcommand{\Lebn}{{\rm Leb}_n} 
\newcommand{\Prob}{{\rm Prob\,}} 
\newcommand{\sinc}{{\rm sinc\,}} 
\newcommand{\ctg}{{\rm ctg\,}} 
\newcommand{\loc}{{\rm loc}} 
\newcommand{\trace}{{\,\,\rm trace\,\,}} 
\newcommand{\Dom}{{\rm Dom}} 
\newcommand{\ifff}{\mbox{\ if and only if\ }} 
\newcommand{\proof}{\noindent {\bf Proof:\ }} 
\newcommand{\remark}{\noindent {\bf Remark:\ }} 
\newcommand{\remarks}{\noindent {\bf Remarks:\ }} 
\newcommand{\note}{\noindent {\bf Note:\ }}

\newcommand{\boldx}{{\bf x}} 
\newcommand{\boldX}{{\bf X}} 
\newcommand{\boldy}{{\bf y}} 
\newcommand{\boldR}{{\bf R}} 
\newcommand{\uux}{\uu{x}} 
\newcommand{\uuY}{\uu{Y}} 
 
\newcommand{\limn}{\lim_{n \rightarrow \infty}} 
\newcommand{\limN}{\lim_{N \rightarrow \infty}} 
\newcommand{\limr}{\lim_{r \rightarrow \infty}} 
\newcommand{\limd}{\lim_{\delta \rightarrow \infty}} 
\newcommand{\limM}{\lim_{M \rightarrow \infty}} 
\newcommand{\limsupn}{\limsup_{n \rightarrow \infty}} 
 
\newcommand{\ra}{ \rightarrow }

\newcommand{\ARROW}[1] 
  {\begin{array}[t]{c}  \longrightarrow \\[-0.2cm] \textstyle{#1} \end{array} } 
 
\newcommand{\AR} 
 {\begin{array}[t]{c} 
  \longrightarrow \\[-0.3cm] 
  \scriptstyle {n\rightarrow \infty} 
  \end{array}} 
 
\newcommand{\pile}[2] 
  {\left( \begin{array}{c}  {#1}\\[-0.2cm] {#2} \end{array} \right) } 
 
\newcommand{\floor}[1]{\left\lfloor #1 \right\rfloor} 
 
\newcommand{\mmbox}[1]{\mbox{\scriptsize{#1}}} 
 
\newcommand{\ffrac}[2] 
  {\left( \frac{#1}{#2} \right)} 
 
\newcommand{\one}{\frac{1}{n}\:} 
\newcommand{\half}{\frac{1}{2}\:} 
 
\def\le{\leq} 
\def\ge{\geq} 
\def\lt{<} 
\def\gt{>} 
 
\def\squarebox#1{\hbox to #1{\hfill\vbox to #1{\vfill}}} 
\newcommand{\qed}{\hspace*{\fill} 
           \vbox{\hrule\hbox{\vrule\squarebox{.667em}\vrule}\hrule}\bigskip} 
 
\title{A Necessary and Sufficient Condition for Invertibility of
  Adapted Perturbations of Identity on Wiener Space\footnote{To appear
  in Comptes Rendus Math\'matiques, Vol. 346}}

\author{ A. S. \"Ust\"unel} 
\date{ } 
\maketitle 
\noindent 
{\bf Abstract:}{\small{ Let $(W,H,\mu)$ be the classical  Wiener space, assume 
that $U=I_W+u$ is an adapted perturbation of identity  satisfying the Girsanov
identity. Then, $U$  is invertible if and only if the kinetic energy of $u$ is
equal to the relative  entropy of the measure induced with the action of $U$
on the Wiener measure $\mu$, in  other words $U$  is invertible if and
only if
$$
\half \int_W|u|_H^2d\mu=\int_W \frac{dU\mu}{d\mu}\log\frac{dU\mu}{d\mu}d\mu\,.
$$
}}\\ 
 
\vspace{0.5cm} 
\noindent 
{\small{\bf Une condition n\'ecessaire et suffisante pour
    l'inversibilit\'e de perturbations d'identit\'e adapt\'ees sur
    l'espace de Wiener}} \\  
 
\noindent 
{\bf Resum\'e:}{\small{ Soit $(W,H,\mu)$ l'espace de Wiener classique, 
et soit $U=I_W+u$ une perturbation d'identit\'e adapt\'ee
satisfaisant \'a l'identit\'e de Girsanov. Alors $U$ est inversible si
et seulement si l'\'energie cin\'etique de $u$ est \'egale \`a
l'entropie relative de la mesure induite par l'action de $U$ sur la
mesure de Wiener $\mu$; en d'autres termes $U$ est inversible si et
seulement si
$$
\half \int_W|u|_H^2d\mu=\int_W \frac{dU\mu}{d\mu}\log\frac{dU\mu}{d\mu}d\mu\,.
$$
}} 
 
\section{Version fran\c{c}aise abr\'eg\'ee} 
Soit $(W,H,\mu)$ l'espace  de Wiener classique: 
$W=C([0,1],\R^d),\,d\geq 1$ $\mu$ est la mesure gaussienne standard et $H$ 
est l'espace de Cameron-Martin dont le produit scalaire et la norme sont 
not\'es respectivement par $(h,k)_H=\int_0^1\dot{h}_s\cdot\dot{k}_sds$
et par  $|\cdot|_H$. Nous noterons par $(\calF_t,t\in [0,1])$ la
filtration canonique du mouvement brownien canonique, eventuellement
compl\'et\'ee. Soit maintenant $u\in L^0_a(\mu,H)$, o\`u cette
derni\`ere repr\'esente les classes d'\'equivalences de variables
al\'eatoires \`a valeurs dans $H$ telles que leurs d\'eriv\'ees
temporelles sont adapt\'ees \`a la filtration brownienne. On notera
par $\rho(u)$ l'exponentielle de Girsanov d\'efinie par
$$
\rho(u)=\exp\left(\int_0^1\dot{u}_s.dW_s-\frac{1}{2}\int_0^1|\dot{u}_s|^2ds\right)
$$
 Nous avons le r\'esultat suivant:
\begin{theorem}
\label{thm-1}
Supposons que $E[\rho(-u)]=1$ et notons par $U$ l'application
$I_W+u$. Alors les propri\'et\'es suivantes sont \'equivalentes:
\begin{enumerate}
\item L'application $U$ est p.s. inversible, son inverse peut
  s'\'ecrire comme $V=I_W+v$ avec $v\in L^0_a(\mu,H)$,
\item L'\'equation differentielle stochastique suivante
\beaa
dV_t&=&-\dot{u}_t\circ V\,dt+dW_t\\
V_0&=&0
\eeaa
poss\`ede une
  solution (forte)  unique,
\item On a l'identit\'e 
$$
\half \int_W|u|_H^2d\mu=\int_W \frac{dU\mu}{d\mu}\log\frac{dU\mu}{d\mu}d\mu\,.
$$
\end{enumerate}
\end{theorem}
L'\'equivalence entre (1) et (2) a d\'ej\`a d\'emontr\'ee dans
\cite{INV} et le reste de cet article sera consacr\'e \`a la preuve de
l'\'equivalence entre (1) et (3).
 
\section{Main results} 
Let $(W,H,\mu)$ the classical Wiener space on $\R^d$, we begin with
the following proposition whose proof follows from the Girsanov theorem
\begin{proposition}
\label{prop-1}
Assume that $L=\rho(- v)$, where $v\in L_a^0(\mu,H)$, i.e.,
$\dot{v}$ is adapted and $\int_0^1|\dot{v}_s|^2ds<\infty$ a.s.
Then there exists $U=I_W+u$, with $u:W\to H$ adapted such that
$U\mu=L\mu$ and $E[\rho(- u)]=1$ if and only if the
following condition is satisfied:
\begin{eqnarray}
\label{MA-1}
1&=&L_t\circ U \,\,E\left[\rho(-u^t)|\calU_t\right]\\
&=&L_t\circ U \,\,E\left[\rho(- u)|\calU_t\right]
\end{eqnarray}
almost surely for any $t\in [0,1]$, where $u^t$ is defined as
$u^t(\tau)=\int_0^{t\wedge \tau}\dot{u}_sds$  and  $\calU_t$ is
the sigma algebra generated by $(w(\tau)+u(\tau),\,\tau\leq t)$.
\end{proposition}
 
\noindent 
Let us calculate $E[\rho(-u^t)|\calU_t]=E[\rho(- u)|\calU_t]$ in terms of the
innovation process associated to $U$. Recall that the term
innovation, which originates from the filtering theory is defined
as (cf.\cite{FKK} and  \cite{BOOK})
$$
Z_t=U_t-\int_0^t E[\dot{u}_s|\calU_s]ds
$$
and it is a $\mu$-Brownian motion with respect to the filtration
$(\calU_t,t\in [0,1])$. A similar proof as the one in \cite{FKK} shows
that any martingale with respect to the filtration of $U$ can be
represented as a stochastic integral with respect to $Z$. Hence, by
the positivity assumption, $E[\rho(- u)|\calU_t]$ can be written
as an exponential martingale
$$
E[\rho(-u)|\calU_t]=\exp\left(-\int_0^t(\dot{\xi}_s,dZ_s)-\frac{1}{2}\int_0^t|\dot{\xi}_s|^2ds\right)\,.
$$
A double utilization of the Girsanov theorem gives the following
explicit result:
\begin{proposition}
\label{prop-2}
We have 
\begin{equation}
\label{cond-exp} E[\rho(-\delta
u)|\calU]=\exp\left(-\int_0^1(E[\dot{u}_s|\calU_s],dZ_s)-\frac{1}{2}\int_0^1|E[\dot{u}_s|\calU_s]|^2ds\right)\,,
\end{equation}
and
\begin{equation}
\label{cond-exp1} E[\rho(-\delta
u)|\calU_t]=\exp\left(-\int_0^t(E[\dot{u}_s|\calU_s],dZ_s)-\frac{1}{2}\int_0^t|E[\dot{u}_s|\calU_s]|^2ds\right)\,,
\end{equation}
almost surely.
\end{proposition}
Combining Propositions \ref{prop-1} and \ref{prop-2}, we obtain
\begin{theorem}
\label{thm-1}
A  necessary and sufficient condition for the relation (\ref{MA-1})
is that
$$
E[\dot{u}_t|\calU_t]=-\dot{v}_t\circ U
$$
$dt\times d\mu$-almost surely.
\end{theorem}
Now we state and prove the main theorem (cf. also \cite{FUZ} for
related problems):
\begin{theorem}
 \label{entropy-thm}
 Assume that $u\in L^2(\mu,H)\cap L^0_a(\mu,H)$ with $E[\rho(-
 u)]=1$. Define $L$ as
$$
L=\frac{dU\mu}{d\mu}=\rho(- v)
$$
where  $v\in L_a^0(\mu,H)$ is given by the It\^o representation
theorem. The map  $U=I_W+u$ is then almost surely
invertible with its inverse  $V=I_W+v$ if and only if
$$
2E[L\log L]=E[|u|_H^2]\,.
$$
In other words, $U$ is invertible if and only if
$$
H(U\mu|\mu)=\half \|u\|^2_{L^2(\mu,H)}\,,
$$
where $H(U\mu|\mu)$ denotes the entropy of $U\mu$ with respect to
$\mu$.
\end{theorem}
\proof
 Since $U$ represents $Ld\mu$, we have, from Theorem \ref{thm-1}, 
$E[\dot{u}_s|\calU_s]+\dot{v}_s\circ U=0$ $ds\times d\mu$-almost
surely. Hence, from the Jensen inequality $E[|v\circ U|_H^2]\leq
E[|u|_H^2]$. Moreover the Girsanov theorem gives
$$
2E[L\log L]=E[|v|_H^2 L]=E[|v\circ
U|_H^2]=E[\int_0^1|E[\dot{u}_s|\calU_s]|^2ds]\,.
$$
Hence the hypothesis implies that
$$
E[|u|_H^2]=E[\int_0^1|E[\dot{u}_s|\calU_s]|^2ds]\,.
$$
From which we deduce that $\dot{u}_s=E[\dot{u}_s|\calU_s]$
$ds\times d\mu$-almost surely. Finally we get
$\dot{u}_s+\dot{v}_s\circ U=0$ $ds\times d\mu$, which is a
necessary and sufficient condition for the invertibility of $U$ (cf.\cite{BOOK,CRAS,INV}) . The necessity is
obvious.

 \qed
\begin{corollary}
With the notations of theorem, $U$ is not invertible if and only if we
have 
$$
\half E[|u|_H^2]>H(U\mu|\mu)\,.
$$
\end{corollary} 
\remark
This result gives an enlightenment about the celebrated counter
example of Tsirelson, cf. \cite{I-W}.

\noindent
{\bf{Acknowledgement:}} Some parts of this work has been done while the
author was visiting the Departement of Mathematics of Bilkent
University, Ankara, Turkey.

\vspace{2cm}
\footnotesize{A. S. \"Ust\"unel, Telecom-Paristech, D\'ept. Infres, \\
46, rue Barrault, 75013, Paris, France\\
ustunel@enst.fr}

\end{document}